\documentclass[11pt]{article}
\begin{document}
\input Vatola.sty
\TagsOnRight

\def\a{\alpha}
\def\b{\beta}
\def\C{{\rm\mbox{l}\!\!\!\mbox{C}}}
\def\Chi{\bar{\bf \xihi}}
\def\chi{\bar{\bf \chi}}
\def\D{{\cal D}}
\def\d{\delta}
\def\e{\varepsilon}
\def\et{\eta}
\def\f{\widetilde{f}}
\def\G{\Gamma}
\def\g{\gamma}
\def\grad{{\rm grad}}

\def\O{\Omega}

\def\intm{\int\limits_M}
\def\into{\int\limits_{\O}}
\def\l{\lambda}
\def\L{\Lambda}
\def\M{\widetilde{M}}
\def\N{\widetilde{N}}
\def\notin{\hspace{2mm}\backslash\hskip-10pt\in}
\def\o{\omega}
\def\p{\partial}
\def\R{{\rm\mbox{I}\!\mbox{R}}}
\def\r{\rho}
\def\s{\sum}
\def\si{\sigma}
\def\t{\tau}
\def\T{{\cal L}}
\def\Tri{\bigtriangledown}
\def\tri{\triangle}
\def\wt{\widetilde}
\def\vf{\varphi}
\def\vf{\varphi}
\def\Z{{\rm\mbox{Z}\!\!\mbox{Z}}}

\begin{center}
{\bf\Large\bf Maps with Prescribed Tension Fields}\\
\vskip 4 mm

{\small Wenyi Chen\ \ \ \ \ \ \  J\"urgen  Jost}

\end{center}


\vskip 8mm

\begin{quote}
{\bf Abstract:}
We consider maps into a Riemannian manifold of nonpositive sectional curvature
with prescribed tension field. We  derive a priori estimates and solve a Dirichlet problem.

   {\bf MR Classification}: 58E20, 53C22.
\end{quote}
\bigbreak
{\Large\bf \S 1.  Introduction}
\bigbreak

Elliptic regularity theory is traditionally concerned with
functions $ f : \Omega \to \R, \Omega$ being some domain in some Euclidean
space $\R^{m}$, or, more generally, in a Riemannian manifold $M$.
$f$ is assumed to solve some elliptic
PDE, and the regularity theory derives estimates of
various norms of $f$ in terms of some data (boundary
conditions, geometry of the domain) and some structural
constants of the elliptic operator. The prototype here is the
Laplace operator, and the elliptic equation in the simplest case then
reduces to the Poisson equation
$$
 \Delta f = v,
$$

\noindent for some prescribed $v$, plus some boundary condition. The
idea of elliptic regularity theory then is to control
some Sobolev or H\"older norm of $f$,
$$
 \Vert f \Vert_{W^{k,p}} , \; \hbox{ or } \Vert f \Vert_{C^{k,\alpha}}
$$
in terms of the corresponding norm of $\Delta f$ for
$k-2 $ in place of $k$ plus some terms depending on
the data. This may then be applied to estimate a
solution of the Poisson equation through the prescribed
right hand side $v$. As  is well-known, such estimates
provide the basis for the existence theory for solutions
of elliptic equations and a guide for the construction of
numerical approximation schemes (for a detailed
presentation, see e.g. \cite{[J2]}).
\medskip

We are interested here not in functions, but in maps
$f : \Omega \to N$ into some Riemannian manifold $N$. This will
make the problem genuinely nonlinear. Nevertheless, the
Laplace operator $\Delta f$ naturally generalizes to the
tension field $\tau (f)$, given in local coordinates by
$$
 \tau(f)^k = \frac{1}{\sqrt{\gamma(x)}} \frac{\partial}{\partial x^{\alpha}}
 \Big( \gamma^{\alpha \beta} (x) \sqrt{\gamma(x)} \frac{\partial f^k (x)}{\partial x^{\beta}} \Big) \\
  + \gamma^{\alpha \beta } (x) \Gamma^i_{jk} (f(x)) \frac{\partial f^i (x)}{ \partial x^{\alpha}}
    \frac{\partial f^j (x)}{\partial x^{\beta}},
$$
with the standard notation:

$(\gamma_{\alpha \beta} (x))_{ \alpha, \beta = 1, \dots, m}$ is the Riemannian metric tensor on the domain
$\Omega \subset M$ w.r.t. local coordinates $x = ( x^1, \dots, x^m),
( \gamma^{\alpha \beta})_{\alpha, \beta = 1, \dots, m} $ its inverse, $ \gamma = \det (\gamma_{\alpha \beta})$.
$\Gamma^i_{jk} (f)$ are the Christoffel symbols of the Riemannian
metric tensor $( g_{ij } (f))_{i,j = 1, \dots, n}$ on $N$, where
we take the liberty to identify $f$ with its local coordinate
representation $(f^1, \dots, f^n)$.
\medskip

Those local coordinates can be chosen in such a manner
that the Christoffel symbols $\Gamma^i_{jk} (f)$ vanish precisely
if $N$ is flat, and in that case, the tension field
reduces to the Laplacian of the domain $\Omega$. In general,
however, a Riemannian manifold is not flat, and so
the tension field then is a nonlinear elliptic operator.
\medskip

In more abstract terms, $\tau (f)$ is a section of the bundle $f^{-1} TN,$
and therefore fundamentally nonlinear as this bundle itself depends on the map $f$.
\medskip

This tension field is the negative gradient field of
the energy functional
$$
 E (f) = \frac{1}{2} \int\vert df(x)\vert^2 \; \hbox{ dvol } (M) \\
 = \frac{1}{2} \int \gamma^{\alpha \beta} (x) g_{ij} (f (x)) \frac{\partial f^i}{\partial x^{\alpha}}
 \frac{\partial f^j}{\partial x^{\beta}} \sqrt{\gamma (x)} dx^1 \dots dx^m
$$
in the same way that the Laplace operator is the negative
gradient field of the Dirichlet integral. Critical points of
the energy integral, i.e. solutions of
$$
 \tau (f) = 0,
$$

\noindent are called harmonic maps and have been intensively
studied.
\medskip

More generally, given some function $\Phi : N \to \R$, one
may look for critical points of the integral
$$
 E_{\Phi} (f) : = \frac{1}{2} \int \vert df (x) \vert^2 \; \hbox{ dvol } (M) + \int \Phi (f(x)) \;
 \hbox{ dvol } (M).
$$
Solutions solve a system of type
$$
 \tau (f) = V (f),
$$
for some vector field $V$ on $N$, and are called
harmonic maps with potential, see e.g. \cite{[CQ]}, \cite{[PW]}, \cite{[FR]}.
\medskip

In the present paper, we shall be concerned
with the system
$$
 \tau (f(x)) = V (x, f(x)),
$$
for some given $V$, without assuming a variational
structure, however. Since $\tau (f)$ is a section of $f^{-1} TN$, also the right
hand side naturally has to depend on $f$, in contrast to the linear Poisson equation
where $V$ is a function of $x$ only.

The existence problem has been studied by
von Wahl \cite{[vW]} in case of a domain $\Omega$ with boundary
and prescribed boundary values, and by Jost-Yau [6]
for this case as well as for the more subtle one of
a closed manifold $M$. In particular, it turns
out that, without a variational structure, on a closed
manifold $M$, the above problem need not possess a
solution. This is in contrast to the
Dirichlet problem on $\Omega$ that has a solution under
general circumstances, as shown in those papers.
\medskip

Actually, all these results need some curvature
assumptions on the target $N$, and, without
imposing a size restriction, one needs the assumption
that $N$ have nonpositive sectional curvature. Therefore,
we shall assume throughout this paper that
{\bf $N$ has nonpositive sectional curvature}. Under
this assumption, we wish to study the above
problem as a nonlinear generalization of the
Poisson equation. Still, while we no longer
have a linear structure, the curvature
assumption implies that there exists an underlying
convex geometry as has been explored in more abstract
terms in \cite{[J1]}.
\medskip

The first part of our paper is concerned with
extending elliptic regularity theory to the present
nonlinear setting. The guiding idea is
that the tension field $\tau$ should assume the role
of the Laplacian $\Delta$. An additional feature, however, is
that not only the geometry is nonlinear, but also that
the topology is nontrivial. Therefore, our estimates
will also involve a topological datum, namely the
homotopy class of the map in question. With that
addition, we are able to generalize the fundamental
elliptic $W^{2,2}$ estimate to our setting. In Theorem. 2.2,
we shall estimate the $L^2$-norms of the first and
second derivatives of any map $f: M \to N$ solely
in terms of the $L^2$-norm of its tension field plus a topological
term that only depends on the homotopy class $[f]$ of $f$.

Of course, the constants in those estimates will depend
on the underlying geometry, in particular on a bound
for the Ricci curvature of $M$.
\medskip

In the second part of our paper, we study the
Dirichlet problem for
$$
 \tau (f) = V.
$$
Our result will need some bound on $V$
depending on the first eigenvalue of $\Omega$. The result
as well as the method are different from those obtained
by von Wahl and Jost-Yau.

\bigbreak

{\Large\bf \S 2.  A nonlinear $W^{2,2}-$ estimate}
\bigbreak

A continuous map $f: M:\rightarrow N$ induces a homeomorphism
$$
\r = f_{\#}: \pi_1(M)\rightarrow\pi_1(N)
$$
of fundamental groups and at the same time a lift
$$
\f: \M \rightarrow\N.
$$
We shall need the $\r-$equivariance of the lift map $\f$,
i.e.
$$
\f(\l x)= \r(\l)\f(x)
$$
for all $x\in \M$, $\l\in\pi_1(M)$ where the fundamental groups $\pi_1(M)$ and $\pi_1(N)$
act isometrically on $\M$ and $\N$ by deck transformations respectively so that $M=\M/\pi_1(M)$
and $N=\N/\pi_1(N)$. There is  a correspondence between the $\r-$equivariant maps from $\M$
 to $\N$ and maps from $M$ to $N$.

    For the complete simply connected Riemannian manifold of nonpositive sectional
 curvature $\N$, the distance function
 $$
d: \N\times \N\rightarrow \R, \ \ d(u,v)=\mbox{the distance from}\;\; u \;\;{\rm to}\;\; v
$$
 is well defined and we have a smooth function $d^2$ on the  manifold
 $\N\times \N$. Let $f_1, f_2:\M\rightarrow\N$ be $\r-$equivariant maps, then the  function
$d^2(f_1, f_2)$ is also smooth on the manifold $M$ because deck transformations are isometric.
For $(u, v)\in \N\times\N$, we choose an orthonormal basis
$e_1, e_2, \cdots, e_n$  for $T_u\N$. By parallel translation along the shortest geodesic from $u$ to $v$
on $\N$, we get a basis $\overline{e}_1,\overline{e}_2, \cdots, \overline{e}_n$  for $T_v\N$. Take
$e_1, e_2, \cdots, e_n, \overline{e}_1,\overline{e}_2, \cdots, \overline{e}_n$ as a local orthonomal frame
for $T_{(u,v)}N\times N$.  Let $\theta^1, \theta^2, \cdots, \theta^m$ be an orthonormal coframe in a neighbourhood
of $x \in \M$, we then  write the differentials of the maps $f_1, f_2: \M\rightarrow\N$ as
 $d f_1 = f^i_{1\a} e_i\otimes\theta^{\a}$, and  $d f_2 = f^{\overline{i}}_{2\a}\overline{e}_i\otimes\theta^{\a}$.
We use the Einstein summation convention.

   Now we define $E(f_1, f_2)$,  the energy of the difference of the two maps $f_1$, $f_2$ by
$$
E(f_1, f_2)=\frac{1}{2}\intm \sum_{i=1}^n \sum_{\a=1}^m\Big(f^i_{1\a}- f^{\overline{i}}_{2\a}\Big)^2 dvol(M)
$$
where the integral is calculated on a fundamental domain of $M$ in $\M$. 
Since 
parallel transport is isometric, we obtain  a triangle inequality
$$
\Big|\sqrt{E(f_1, f_3)}-\sqrt{E(f_3, f_2)}\Big|\leq \sqrt{E(f_1, f_2)},
$$
Inparticular, the energy of the difference is symmetric, and for the
energy of a map $f$,
$$
E(f)=\frac{1}{2}{\intm} \sum_{i=1}^n \sum_{\a=1}^m\Big(f^i_{\a}\Big)^2 dvol(M),
$$
 we then have
$$
\Big|\sqrt{E(f_1)}-\sqrt{E(f_2)}\Big|\leq \sqrt{E(f_1, f_2)}.
\eqno(2.0)
$$

Putting $X_{\a}=f^i_{1\a} e_i+ f^{\overline{i}}_{2\a}\overline{e}_i$, we then are able to express the Laplacian of
$d^2(f_1, f_2)$ as
$$
\Delta d^2(f_1, f_2) = (d^2)_{X_{\a}X_{\a}}+d\{d_if^i_{1\a\a}+d_{\overline{i}}f^{\overline{i}}_{2\a\a}\}
$$
where $\t(f_1)^1=f^i_{1\a\a}$ and $\t(f_1)^2=f^{\overline{i}}_{2\a\a}$ are the components of the tension fields
of the maps $f_1$, $ f_2$ respectively. It was shown in \cite{[SY]} that if $K_N\leq 0$ then
$$
 (d^2)_{X_{\a}X_{\a}}\geq 2\sum_i \Big(f^i_{1\a}- f^{\overline{i}}_{2\a}\Big)^2.
$$
Hence
$$
\Delta d^2(f_1, f_2)\geq 2\sum_i\Big(f^i_{1\a}- f^{\overline{i}}_{2\a}\Big)^2- d\{|\t(f_1)|+|\t(f_2)|\}.
\eqno(2.1)
$$
We have also
$$
\Delta d(f_1, f_2)\geq - \{|\t(f_1)|+|\t(f_2)|\}.
\eqno(2.2)
$$
Let $f_1$ and $f_2$ be maps from a domain $\O$ of the manifold $M$, smooth enough in $\O$ and on the boundary
$\p\O$. If two maps $f_1$ and $f_2$ coincide on the boundary $\p\O$, we then have by (2.2)
$$
\into|\nabla d(f_1, f_2)|^2 d vol\leq \Big\|d(f_1, f_2)\Big\|_{L^2(\O)}
\Big\{\Big\|\t(f_1)\Big\|_{L^2(\O)}+\Big\|\t(f_2)\Big\|_{L^2(\O)}\Big\}.
$$
Therefore
$$
\Big\|d(f_1, f_2)\Big\|_{L^2(\O)}\leq \l(\O)^{-1}
\Big\{\Big\|\t(f_1)\Big\|_{L^2(\O)}+\Big\|\t(f_2)\Big\|_{L^2(\O)}\Big\}
\eqno(2.3)
$$
 where the constant $\l(\O)$ is the first eigenvalue of the Laplacian on $\O$,  see \cite{[DW]}.

On the other hand, the integral of (2.1) over $\O$ gives
$$
E(f_1, f_2)\leq 2^{-1} \into\Delta d^2(f_1, f_2) + d(f_1, f_2)\{|\t(f_1)|+|\t(f_2)|\}.
$$
Notice that $d^2(f_1, f_2)\geq 0$ in $\O$ and $d^2(f_1, f_2)=0$ on $\p\O$, the derivative of
$d^2(f_1, f_2)$ along the outer normal vector of $\O$ at $\p\O$ will be nonpositive. Hence
the Stokes formula gives
$$
\into\Delta d^2(f_1, f_2)\leq 0.
$$
We then obtain an energy estimate for two maps that have the same  boundary values and belong to the same homotopy
class,  with (2.3)
$$
E(f_1, f_2)\leq  \l(\O)^{-1}
\Big\{\Big\|\t(f_1)\Big\|^2_{L^2(\O)}+\Big\|\t(f_2)\Big\|^2_{L^2(\O)}\Big\}.
\eqno(2.4)
$$
For the case of a compact manifold $M$ without boundary, we need a deeper analysis.
 Integrating the inequality (2.1) over a fundamental domain, we get
$$
E(f_1, f_2)\leq \frac{1}{4} \Big\|d(f_1, f_2)\Big\|_{L^2(M)}
\Big\{\Big\|\t(f_1)\Big\|_{L^2(M)}+\Big\|\t(f_2)\Big\|_{L^2(M)}\Big\}.
\eqno(2.5)
$$

A direct corollary of the triangle inequality and the inequality (2.2) is that if $h_1, h_2$ are two
homotopic harmonic maps, then $E(h_1, h_2)=0$ and $E(f, h_1)=E(f, h_2)$ whenever $[f]=[h_1]$.

{\bf Lemma 2.1}:
{\it Let $f_0$, $f_1: \M\rightarrow \N$ be $\r-$equivariant maps. Define $f_t:\M\rightarrow \N$ by exponential map
$$
f_t(x)=\exp_{f_0(x)}\{t\exp^{-1}_{f_0(x)} f_1(x)\},
$$
then $f_t$ is also $\r-$equivariant and $\sqrt{E(f_t)}$
 is a convex function of $t$.
}

{\bf Proof}:
 Let $\g_x$ be the geodesic on $\N$ from $f_0(x)$ to  $f_1(x)$, then $\g_x$ is also  $\r-$equivariant.
Hence, $f_t$ is  $\r-$equivariant. Fix $\a\in \{1, 2, \cdots, n\}$, we claim that $\nabla_{\a} f_t(x)$ is
a Jacobi field along the geodesic $\g_x$. In fact, let $\nu_{\a}\in T_x \M$ be the vector dual to
$\theta^{\a}$, $c(s)$ be the geodesic on $\M$ with $c(0)=x$, $c^\prime(0)=\nu_{\a}$,
$C:(0, 1)\times(-\delta, \delta)\rightarrow \N$, $C(t, s)= \exp_{f_0(c(s))}
\{t \exp^{-1}_{f_0(c(s))} f_1(c(s))\}$, then
for any fixed $s$, $C(\cdot, s)$ is a geodesic. On the other hand
$$
\nabla_{\a} f_t(x)=\frac{\textstyle \p}{\textstyle \p s} C(t, s)\Big|_{s=0}.
$$
By the Jacobi field equation, we have
$$
\align
&\frac{\textstyle \p^2}{\textstyle \p t^2}\sqrt{E(f_t)} = \frac{\textstyle -1}{\textstyle \sqrt{E(f_t)}}
\intm \langle R(\nabla_{\a} f_t(x), \stackrel{\cdot}{\g})\stackrel{\cdot}{\g}, \nabla_{\a} f_t(x)\rangle d vol(M)\\
&\;\; +\frac{\textstyle 1}{\textstyle \sqrt{E(f_t)}}
\intm\sum_{\a}\Big|\frac{\textstyle \p}{\textstyle \p t}\nabla_{\a} f_t(x)\Big|^2 dvol (M)\\
&\;\;-\frac{\textstyle 1}{\textstyle E(f_t)^{\frac{3}{2}}}
\Big|\intm\nabla_{\a} f_t(x)\nabla_{\a}\frac{\textstyle \p}{\textstyle \p t} f_t(x)d vol(M)\Big|^2 \\
& \;\; \geq 0.
\endalign
$$
 As required.
\bigbreak

{\bf Lemma 2.2}:
{\it
There exists a constant $C$  which depends only on the homotopy class $[f]$ such that
$$
\Big\|d(f, h)\Big\|^2_{L^2(M)}\leq C E(f, h)
\eqno(2.6)
$$
holds for a harmonic map in the homotopy class $[f]$.
}

{\bf Proof} We prove the  Lemma by contradiction. If  there were a homotopy class such that
the inequality (2.3) did not hold, we then could
find a sequence $f_k$ of maps and harmonic maps $h_k$ in the same homotopy class such that
$$
\Big\|d(f_k, h_k)\Big\|^2_{L^2(M)}=\inf_{\t(h)=0}\Big\|d(f_k, h)\Big\|^2_{L^2(M)}\geq k^2  E(f_k, h_k)
\eqno(2.7)
$$
with $\Big\|d(f_k, h)\Big\|^2_{L^2(M)}\geq 1 $. Notice that the energy of the difference $E(f_k, h)$ is
 independent of the choice of the harmonic map $h$. Define
$$
f_k^t(x)=\exp_{h_k(x)}\{t\exp^{-1}_{h_k(x)} f_k(x)\}.
$$
We have
$$
\align
&\Big\|d(f^t_k, h_k)\Big\|^2_{L^2(M)}=t^2\Big\|d(f_k, h_k)\Big\|^2_{L^2(M)}\geq k^2 t^2 E(f_k, h_k)\\
&\;\; \geq k^2 t^2 \Big(\sqrt{E(f_k)}-\sqrt{E(h_k)}\Big)^2.
\endalign
$$
It follows from the convexity  of the energy $E(f^t_k)$ that
$$
\Big\|d(f^t_k, h_k)\Big\|_{L^2(M)}
\geq k \Big|\sqrt{E(f_k)}-\sqrt{E(h_k)}\Big|.
\eqno(2.8)
$$
Choose $t=t_k$ such that $\Big\|d(f^t_k, h_k)\Big\|_{L^2(M)}=1$. From (2.8), $f^{t_k}_k$ is a minimizing sequence for the energy.
Therefore the sequence $f^{t_k}_k$ converges strongly to a harmonic map $h$, i.e.
$$
\Big\|d(f^{t_k}_k, h)\Big\|_{L^2(M)} \rightarrow 0.
\eqno(2.9)
$$
Notice that
$$
\align
&\Big\|d(f_k, f^{t_k}_k)\Big\|_{L^2(M)}=\Big\|d(f_k, h_k)-d(h_k, f^{t_k}_k)\Big\|_{L^2(M)}\\
&=\Big\|(1-t_k)d(f_k, h_k)\Big\|_{L^2(M)}\\
&=\Big\|d(f_k, h_k)\Big\|_{L^2(M)}-t_k\Big\|d(f_k, h_k)\Big\|_{L^2(M)}\\
&=\Big\|d(f_k, h_k)\Big\|_{L^2(M)}-1.
\endalign
$$
Hence for $k$ large enough,
$$
\align
&\Big\|d(f_k, h)\Big\|_{L^2(M)}\\
&\leq \Big\|d(f_k, f^{t_k}_k)\Big\|_{L^2(M)}+\frac{1}{2} \\
&=\Big\|d(f_k, h_k)\Big\|_{L^2(M)}-\frac{1}{2}.
\endalign
$$
This is a contradiction to the choice of the harmonic maps $h_k$. This proves that the inequality (2.3)
holds for  $\inf_{\t(h)=0}\Big\|d(f, h)\Big\|^2_{L^2(M)}\geq 1$. Notice that the harmonic map and the inequality
(2.3) are invariant under the  rescaling of the metric on $N$, the inequality (2.3) holds for all maps in the same
homotopy class.\\\\

Now we summarize the above, that is, (2.0), (2.5) and (2.6) as:

{\bf Theorem 2.1}:
{\it
 Let $M$  be a compact Riemanian manifold with or without boundary, $N$  a compact Riemannian manifold
 of nonpositive sectional curvature. For a given homotopy class of  maps from $M$ to $N$, with prescibed boundary values when the boundary of $M$ is not empty, there is a constant
  $C$ such that for any map $f$ in this homotopy class there is a harmonic map $h$ in
the homotopy class such that
  $$
  \Big(\intm |d f|^2\Big)^{\frac{1}{2}} \leq \Big(\intm |d h|^2\Big)^{\frac{1}{2}}
 + C\Big(\intm |\t(f)|^2\Big)^{\frac{1}{2}}
  \eqno(2.10)
$$
  where
  $\t(f)$ is the tension field.

}
The first term on the right hand side  depends only on the homotopy class of $f$, because all harmonic maps
in the same homotopy class have the same -minimal- energy.

Recall the Bochner formula (\cite{[JJ]}, Ch. 8)
$$
\align
&\Delta e(f)(x) = |\nabla d f|^2+\langle\d(\tau(f)), df\rangle  +
\langle df(Ric^M(e_{\a})), df(e_{\a})\rangle \\ & \ \ \ -\langle R^N(df(e_{\a}), df(e_{\beta}))df(e_{\beta}),
 df(e_{\a})\rangle .
\endalign
$$
Integrating it on the domain $M$ and making use of the inequality (2.10) under the
assumption of nonpositive curvature, we get the main result of this section, namely an estimate
for the $L^2-$ norm of the first and second derivatives of a map in a given homotopy
class in terms of its tension field.

{\bf Theorem 2.2:}
{\it
Let $M$  be a compact Riemanian manifold  without boundary, $N$  a compact Riemannian manifold
 of nonpositive sectional curvature. For a given homotopy class of maps from $M$ to $N$,
  let $h$ be a harmonic map  in
that homotopy class. We then have, for any map $f$ in that homotopy class, 
$$
\intm |df|^2 dvol(M) +\intm |\nabla d f|^2 d vol(M)\leq C_1\intm |\t(f)|^2 dvol(M)+
  C_2\intm|d h|^2 dvol(M)
$$
where the constant $C_1$ depends only on the homotopy class of the map $f$
and the constant $C_2$ is $ 1+\|Ric^M\|_{\infty}$, where $Ric^M$ is the
Ricci curvature of $M$.
}

We point out once more that the last term on this inequality depends only on the homotopy class (and on the geometry of $M$ and $N$), but not in any way on the map $f$ in that homotopy class that we are estimating here.

By differentiation, we may then also obtain higher order estimates in a standard manner.
\bigbreak

\bigbreak

{\Large\bf \S 3. Boundary value problems}
\bigbreak

We now let $\O$ be a domain in a manifold $M$  with a
nonempty boundary $\p\O$ and compact closure $\overline\O$. For the moment, we assume that the
map $g:\overline\O \rightarrow N$ is of class $C^{2, \a}$.

Consider the parabolic system
$$
\left\{
  \begin{array}{l}
  \tau(f)-\frac{\textstyle \p f}{\textstyle \p t}= V(f),\\
  f(x, 0)=g(x), \ \  x\in \O,\\
  f(x, t)=g(x),  \ \  x\in \p \O,\\
  f: \overline\O\times [0, \infty) \rightarrow N.
 \end{array}
 \right.
 \eqno(3.1)
$$

{\bf Remark:}
If there is an underlying variational structure, as for harmonic maps with potential, namely,
$$
E_{\varphi}(f)=\frac{1}{2}\into |df|^2 dvol(M)+ \into \varphi(f(x)) dvol(M),
$$
for some function $\varphi: N\rightarrow R$, we have the Euler-Lagrange equation as
$$
\tau(f)(x)=\nabla \varphi(f(x)).
$$
We have
$$
\frac{\p }{\p t} E_{\varphi}(f(\cdot, t))= -\into\Big|\frac{\p }{\p t}f(x, t)\Big|^2 dvol(M),
$$
 for a solution of the associated parabolic problem
$$
\frac{\p f}{\p t} = \tau(f)- \nabla \varphi(f),
$$
and therefore, if we assume that $ E_{\varphi}(f)$ is bounded from below, i.e.,
$$
E_{\varphi}(f)\geq C
$$
for some $C\in \R$ and all $f$, then $\frac{\p }{\p t}f(\cdot, t)$ subconverges to zero in $L^2$ for
$t\rightarrow \infty$, and the analysis becomes easy. This has been explained in the literature, see
e. g., \cite{[FR]}. Without a variational structure, however, the problem is more difficult.
\bigbreak
 The smoothness of the map $g$ and the theory of linear parabolic systems give us a short time solution
 of the parabolic problem (3.1), i.e., there is a positive $T$ so that there is a $C^{2, \a}$
 solution on $\O\times [0, T)$. The condition we impose on the vector field is that
 $$
 \begin{array}{l}
 \nabla  V(X, X) \geq -\mu|X|^2, \ \ \ X\in \Gamma(T\O).
\end{array}
 \eqno(3.2)
 $$
We also let $\l(\O)$ be the first nontrivial eignvalue for the Dirichlet problem on $\O$.

{\bf Lemma 3.1:}\ \ \
{\it  Let $f$ satisfy the parabolic system (3.1) where the vector field $V$ satisfies (3.2) with
$\mu\leq \frac{\textstyle 3}{\textstyle 4}\l(\O)$, then
$$
\begin{array}{l}
\into |\frac{\p }{\p t}f(\cdot, t)|^4\leq \into |\tau(g)-V(g)|^4.
\end{array}
\eqno(3.3)
$$
}

{\bf Proof:}\ \ \
A direct computation gives that
$$
\begin{array}{l}
(\Delta - \frac{\p}{\p t})\langle\frac{\p f}{\p t},\frac{\p f}{\p t}\rangle
   =2\nabla V(\frac{\p f}{\p t}, \frac{\p f}{\p t})
   +2\langle\nabla_{\nu_i}\frac{\p f}{\p t},\nabla_{\nu_i}\frac{\p f}{\p t}\rangle \\
   \hskip 6mm- 2\langle R^N(f_{\ast} \nu_i,\frac{\p f}{\p t})f_{\ast} \nu_i,\frac{\p f}{\p t}\rangle .
   \end{array}
\eqno(3.4)
 $$
With the assumptions of the Lemma and the nonpositivity of the sectional curvature of the manifold $N$, it follows from
the above identity
$$
\begin{array}{l}
(\Delta - \frac{\p}{\p t})\Big|\frac{\p f}{\p t}\Big|^2 \geq
   -\frac{\textstyle 3}{\textstyle 2}\l(\O)\Big|\frac{\p f}{\p t}\Big|^2
   +2\Big|\nabla\frac{\p f}{\p t}\Big|^2.
  \end{array}
\eqno(3.5)
 $$
Notice that
$$
2\Big|\nabla\frac{\p f}{\p t}\Big|^2\cdot\Big|\frac{\p f}{\p t}\Big|^2
\geq \frac{1}{2}\Big|\nabla|\frac{\p f}{\p t}|^2\Big|^2.
$$
Multiplying the  two sides of the inequality (3.4) by $\Big|\frac{\p f}{\p t}\Big|^2$
and then integrating over $\O$, one obtains
$$
-\frac{1}{2}\frac{\p}{\p t}\into\Big|\frac{\p f}{\p t}\Big|^4
\geq \frac{3}{2}\into\Big|\nabla|\frac{\p f}{\p t}|^2\Big|^2-
\frac{\textstyle 3}{\textstyle 2}\l(\O)\into\Big|\frac{\p f}{\p t}\Big|^4.
$$
Because $\frac{\p f}{\p t}=0$ on the boundary of the domain $\O$, the right side of the above
inequality is nonnegative. Hence
$$
\frac{\p}{\p t}\into\Big|\frac{\p f}{\p t}\Big|^4\leq 0.
$$
This proves Lemma 3.1.

Let us return to the inequality (3.5)
$$
(\Delta - \frac{\p}{\p t})\langle\frac{\p f}{\p t},\frac{\p f}{\p t}\rangle
 \geq -C\langle\frac{\p f}{\p t},\frac{\p f}{\p t}\rangle .
$$
Set $\psi(x, t)=\exp\{Ct\}\langle\frac{\p f}{\p t},\frac{\p f}{\p t}\rangle $, then
$$
\begin{array}{l}
(\Delta - \frac{\p}{\p t})\psi\geq 0.
\end{array}
\eqno(3.6)
$$
and $\psi = 0$ on the boundary of $\O$.
By Moser iteration (cf \cite{[JM]}), one obtains that
$$
\psi^2(x, t_1)\leq C \Big(1+\frac{1}{t_1-t_0}\Big)\int_{\O\times (t_0, t_1)} \psi^2(x, t)
$$
with the constant $C$ independent of $t$. It follows from Lemma 3.1 and the above inequality that
$$
\Big|\frac{\p f}{\p t}\Big|^4\leq C\into |\tau(g)-V(g)|^4.
$$
Therefore we get a uniform bound for $|\frac{\p f}{\p t}|$. Now we have
$$
\Delta d(f(\cdot,t), g(\cdot))\geq  -C
\eqno(3.7)
$$
where the constant $C$ comes from the uniform bound of $|\frac{\p f}{\p t}|$ and the vector
field $V$ by (2.2) and the parabolic system (3.1). A direct corollary of the inequality (3.7)
is that the distance function $d(f(\cdot,t), g(\cdot))$ has a uniform bound which is
independent of $t$. Indeed applying the
weak maximum principle to the elliptic inequality (3.7) one gets
$$
\sup_{\O}d(f(\cdot,t), g(\cdot))\leq \sup_{\p\O}d(f(\cdot,t), g(\cdot))+ C|\O|^{\frac{1}{m}}.
$$

Consider the Dirichlet problem for the inhomogeneous Laplace equation in $\O$
$$
\Big\{
\begin{array}{ll}
\Delta u = -C,\ \ \ \ &\mbox{in} \ \ \O;\\
u = 0,\ \ \ \ \ \ &\mbox{on} \ \ \p\O ,
\end{array}
\Big.
\eqno(3.8)
$$
where the constant  comes from (3.7).

Let $u$ be the solution of problem (3.8), then
$$
\Delta\{ d(f(\cdot,t), g(\cdot)) - u(\cdot)\} \geq 0.
$$
Hence
$$
\sup_{\O}\{d(f(\cdot,t), g(\cdot)) - u(\cdot)\}\leq \sup_{\p\O}\{d(f(\cdot,t), g(\cdot)) - u(\cdot)\}.
$$
Therefore
$$
d(f(\cdot,t), g(\cdot)) \leq  u(\cdot).
\eqno(3.9)
$$

The solution of the problem (3.8) is of class $C^{2,\a}$ if the boundary $\p\O$ smooth enough. So we have
$$
u(x) \leq  C d(x, \p\O).
$$
This implies
$$
d(f(\cdot,t), g(\cdot)) \leq  C d(x, \p\O).
\eqno(3.10)
$$
Take  $\nu_1, \cdots, \nu_m$ to be a local orthonomal frame at $z_0\in\p\O$ with $\nu_1, \cdots, \nu_{m-1}$
 the tangent vectors of $\p\O$ and  $\nu_m$  the normal direction of the boundary. At the point $z_0\in\p\O$
we  then have
$$
\frac{\p f(\cdot,t)}{\p \nu_j}=\frac{\p g(\cdot)}{\p \nu_j}, \ \ \ 1\leq j\leq m-1,
$$
by the boundary condition, and
$$
\Big|\frac{\p f(\cdot,t)}{\p \nu_m}\Big|\leq C
$$
because of (3.10).
Therefore we obtain a uniform bound for the gradient of the maps $f(\cdot,t)$, i.e.,
$$
e(f)(x)\leq C, \ \ \  x\in\p\O
$$
where the constant $C$ is independent of $t$.
\vskip 0.5cm
>From the above argument we know there is a constant $R$ so that
$$
d(f(\cdot,t), g(\cdot)) \leq  R \ \ \ \  \mbox{on}\ \O.
$$
Let $F= 2R^2-d^2(f(\cdot,t), g(\cdot))$, then $F\geq R^2$. Set
$$
A(x, t)=\frac{e(f)(x,t)}{F^2}.
$$
We have
$$
(\Delta - \frac{\p}{\p t})A(x, t)=\frac{(\Delta - \frac{\p}{\p t})e(f)}{F^2}
-\frac{2e(f)(\Delta - \frac{\p}{\p t})F}{F^3}-\frac{4\nabla e(f)\cdot\nabla F}{F^3}
+\frac{6e(f)|\nabla F|^2}{F^4}.
\eqno(3.11)
$$
By the Bochner formula, (cf \cite{[JJ]}, Ch. 8),
$$
\begin{array}{l}
(\Delta - \frac{\p}{\p t}) e(f) = \langle (\Delta - \frac{\p}{\p t})\frac{\p}{\p x_i}f, \frac{\p}{\p x_i}f\rangle \\
 \hskip 3mm =\langle\nabla_{f_{\ast}{\nu_i}} V, f_{\ast}{\nu_i}\rangle
   +\langle\nabla_{\nu_i\nu_j} f, \nabla_{\nu_i\nu_j}f\rangle \\
   \hskip 6mm- \langle R^N(f_{\ast} \nu_i,f_{\ast} \nu_j)f_{\ast} \nu_i,f_{\ast} \nu_j\rangle
     +\langle f_{\ast}Ric^M \nu_i, f_{\ast} \nu_i\rangle \\
   \hskip 3mm
   \geq |B(f)|^2-C e(f),
 \end{array}
 \eqno(3.12)
 $$
 where $B(f)=\nabla_{\nu_i\nu_j} f$ is the Hessian of the maps $f(\cdot, t)$.
On the other hand,
$$
\begin{array}{l}
(\Delta - \frac{\p}{\p t})F= \nabla^2_N F(f_{\ast}\nu_i, f_{\ast}\nu_i)-\langle\nabla_N F, \tau(f)-\frac{\p f}{\p t}\rangle \\
\hskip 1cm =\nabla^2_N F(f_{\ast}\nu_i, f_{\ast}\nu_i)-\langle\nabla_N F, V\rangle .
\end{array}
$$
By  Jacobi field estimates, (see, e.g. \cite{[JJ]}, Ch. 4)
$$
\nabla^2_N F(f_{\ast}\nu_i, f_{\ast}\nu_i)\leq -2 e(f).
$$
 Hence
$$
(\Delta - \frac{\p}{\p t})F\leq -2 e(f)+2 R |V|.
$$
Here we have used the property that $|\nabla_N F|=2d(f(\cdot,t), g(\cdot))\leq 2 R$.
Returning to (3.11), we have
$$
(\Delta - \frac{\p}{\p t})A(x, t)\geq \frac{4}{F^3}e^2(f)
-\frac{C}{F^2}e(f)-\frac{4R|V|}{F^3} e(f)+I,
\eqno(3.13)
$$
where
$$
I=\frac{1}{F^2}|B(f)|^2-\frac{4}{F^3}\nabla e(f)\cdot \nabla F+ \frac{6e(f)}{F^4}|\nabla F|^2.
$$
Because
$$
\nabla A=\frac{F\nabla e-2 e \nabla F}{F^3},
$$
we have
$$
I=\frac{1}{F^2}|B(f)|^2-\frac{2}{F^3}|\nabla e(f)\cdot \nabla F+ \frac{2e(f)}{F^4}|\nabla F|^2
-\frac{2\nabla A\cdot \nabla F}{F}.
$$
Notice that
$$
|B(f)|^2-\frac{2}{F}|\nabla e(f)\cdot \nabla F+ \frac{2e(f)}{F^2}|\nabla F|^2
\geq \Big(|B(f)|-|df||\nabla F|\Big)^2,
$$
we have
$$
I\geq-\frac{2\nabla A\cdot \nabla F}{F}.
$$
Therefore
$$
\left(\Delta - \frac{\p}{\p t}\right)A(x, t)\geq 
\frac{4}{F^3}e^2(f)-\frac{C}{F^2}e(f)-
\frac{4R|V|}{F^3} e(f)-\frac{2\nabla A\cdot \nabla F}{F}.
$$
Or
$$
(\Delta - \frac{\p}{\p t})A(x, t)\geq 4F A^2
-CA-\frac{4R|V|}{F}A-\frac{2\nabla A\cdot \nabla F}{F}.
\eqno(3.14)
$$
Let
$$
A(t)=\max_{x\in\overline{\O}} A(x, t),
$$
and assume that $A(t)=A(x_t, t)$ for some point $x_t$. If $x_t\in\p\O$ for some $t$, we have
$$
e(f)(x, t)\leq 2R^2 A(x, t)\leq 2R^2 A(x_t, t)\leq 2 e(f)(x_t, t).
$$
It follows from (3.11) that
$$
e(f)(x, t)\leq 2C
\eqno(3.15)
$$
where the constant $C$ is the same as in (3.11). On the other hand if $x_t\in \O$, we have
$$
 - \frac{\p}{\p t}A|_{(x_t, t)}\geq 4F A^2-CA-\frac{4R|V|}{F}A
$$
Without loss of generality, we may assume that $\frac{\p}{\p t}A|_{(x_t, t)}\geq 0$. Hence
$$
4F A^2(t)-CA(t)-\frac{4R|V|}{F}A(t)\leq 0.
$$
That is to say
$$
A(x, t)\leq \frac{C}{4F}+\frac{R|V|}{F^2}.
$$
Therefore
$$
e(f)(x,t)\leq CR^2+R\|V\|,
\eqno(3.16)
$$
where $\|V\|=\max_{y\in N} |V(y)|$.

Up to now, we have proved

{\bf Theorem 3.1:}
{\it Let the vector field $V$ satisfy the condition (3.2) with
$\mu\leq \frac{3}{4}\l(\O)$ for the first eigenvalue $\l(\O)$
of the Laplacian in $\O$. If the maps $f(\cdot, t)$ satisfy
the parabolic system (3.1) for $0\leq t < T$, then there is
a constant $C$ which is independent of $t$ so that
$$
|df(\cdot, t)|,\ \ |\frac{\p f}{\p t}|\leq C.
\eqno(3.17)
$$
}

By linearizing and using the theory of linear  parabolic systems and the implicit function theorem, one gets from (3.17) that
(3.1) has a solution for all of $t$. That is

{\bf Corollary 3.1:}
{\it Let the vector field $V$ satisfy the condition (3.2) with
$\mu\leq \frac{3}{4}\l(\O)$ for the first eigenvalue $\l(\O)$
of the Laplacian in $\O$, then
the parabolic system (3.1) has a solution for all time $t\in [0, \infty)$. Moreover, the solution of
the parabolic system (3.1) has a uniform $C^{2, \a}$ bound.
}

The last claim of Corollary 3.1 means that   any sequence  $t_k\rightarrow \infty$ will contain a subsequence
$t_{n^{\prime}}\rightarrow \infty$ so that $f(\cdot, t_{k^{\prime}})$  converges to a map $f$ in $C^{2}$. This leads
to the existence for the Dirichlet problem
$$
\left\{
  \begin{array}{l}
  \tau(f)= V(f),\\
  f(x)=g(x),  \ \  x\in \p \O\\
  \end{array}
 \right.
 \eqno(3.18)
$$
under a somewhat stronger condition than Theorem 3.1.

{\bf Theorem 3.2:}
{\it Let the vector field $V$ satisfy the condition (3.2) with
$\mu\leq \frac{3}{4}\l(\O)-\e$ for the first eigenvalue $\l(\O)$
of the Laplacian in $\O$, then the Dirichlet problem (3.18) has a
solution.
}

{\bf Proof:}
By the same calculation  as in the proof of Lemma 3.1, one gets that
 $$
-\frac{1}{2}\frac{\p}{\p t}\into\Big|\frac{\p f}{\p t}\Big|^4
\geq \frac{3}{2}\into\Big|\nabla|\frac{\p f}{\p t}|^2\Big|^2-
\Big(\frac{\textstyle 3}{\textstyle 2}\l(\O)-2\e\Big)\into\Big|\frac{\p f}{\p t}\Big|^4.
$$
The Poincar\'e inequality gives
$$
\frac{\p}{\p t}\into\Big|\frac{\p f}{\p t}\Big|^4\leq -4\e\into\Big|\frac{\p f}{\p t}\Big|^4.
$$
That is to say
$$
\into\Big|\frac{\p f}{\p t}\Big|^4\leq C_0 \exp\{-4\e t\}.
$$
Hence, for a sequence  $t_k$ with $f(\cdot, t_k)$ convergent in $C^{2}$,
$$
\lim_{t_k\rightarrow \infty}\frac{\p f}{\p t}(x,t_k)=0.
$$
Therefore the limit $f$ of $f(\cdot, t_k)$ will solve the problem (3.18).

{\bf Acknowledgement:}  The first author would like to express his gratitude
to the Max Planck Institute for Mathematics in the Sciences for support and excellent
working conditions.

\vskip 4mm


Wenyi Chen

{\small

Department of Mathematics

Wuhan University

Wuhan P. R. China 430072
}

\vskip 4mm

J\"urgen  Jost

{\small

Max Planck Institute for Mathematics in the Sciences.

Inselstr. 22$-$26,

D$-$04103 Leipzig

Germany
}

\end{document}